\documentclass{article}

\usepackage{amssymb}
\usepackage{amsmath}

\usepackage[english,francais]{babel}

\newtheorem{theorem}{Theorem}[section]

\newtheorem{e-proposition}[theorem]{Proposition}
\newtheorem{corollary}[theorem]{Corollary}
\newtheorem{e-definition}[theorem]{Definition\rm}

\newcommand{\N}{{\mathbb{N}}}   
\newcommand{\Z}{{\mathbb{Z}}}   
\newcommand{\Q}{{\mathbb{Q}}}   
\newcommand{\U}{{\mathcal{U}}}   
\newcommand{\D}{{\mathcal{D}}}   
\newcommand{\vp}{\vspace{0,5cm}}   
\newcommand{\G}{{\mathcal{G}}}

\newcommand{\h}{{\mathbf{H}}}
\newcommand{\q}{{\mathbb{L}}}
\newcommand{\f}{{\mathbb{F}}}
\newcommand{\hp}{{\mathbb{H}}}    
\newcommand{\n}{{\eta}}
\def\og{\leavevmode\raise.3ex\hbox{$\scriptscriptstyle\langle\!\langle$~}}
\def\fg{\leavevmode\raise.3ex\hbox{~$\!\scriptscriptstyle\,\rangle\!\rangle$}}

\begin{document}
\centerline{}


\selectlanguage{english}
\begin{center}
\Large \textbf{Hyperbolicity of orders of quaternion algebras and group rings}
\normalsize $$\textrm{S. O. Juriaans}^{a} \qquad \textrm{A. C. Souza
  Filho}^{b}$$
\end{center}


\selectlanguage{english}

\begin{center}
{\it Instituto de Matem\'atica e Estat\'\i stica,
 Universidade de S\~ao Paulo,
Caixa Postal 66281, S\~ao Paulo, CEP 05315-970 - Brazil}
\end{center}
\begin{center}
\small
email adresses:
$^{a}$ostanley@ime.usp.br $\quad$ $^{b}$calixto@ime.usp.br
\end{center}
\normalsize



\begin{abstract}
\selectlanguage{english}
For a given division algebra of the quaternions, we construct two types of units of its
$\Z$-orders: Pell units and Gauss units. Also, if $K=\Q \sqrt{-d}$ for $d$ a
square free and $R=I_K$, we classify $R$ and $G$ such that $\U_1(RG)$ is
hyperbolic. In particular, with a suitable geometric approach we prove that $\U_1(RK_8)$ is
hyperbolic iff $d>0$ and $d \equiv 7 \pmod 8$. In this case, the hyperbolic
boundary $\partial(U_1(RG))\cong S^2$, the two dimensional sphere.
\end{abstract}
\vskip 0.5\baselineskip

\selectlanguage{francais}


\selectlanguage{francais}

\selectlanguage{english}


\section{Introduction}
\label{}

Hyperbolic groups were defined firstly by Gromov \cite{grmv}, from the concept
of hyperbolic metric space. 

\vp

Let $G$ be a finitely generated group and $\G$ its Cayley graph
 with the lenght metric, $G$ is hyperbolic if $\G$ is hyperbolic.

\vp

Gromov showed that if $\Gamma$ is hyperbolic, then it does not contain
a free abelian group of rank two, i.e.,  $\Z^{2} {\not \hookrightarrow}
\Gamma$. If $G$ is finite then $\Q G$ has at most one Wedderburn component that is not a
division ring and it is isomorphic to $M_2 (\Q)$. This was first proved by Jespers in
\cite{jsp}. Still in \cite{jsp}, Jespers classifies the finite groups $G$ with non abelian
free normal complement in  $\U_1(\Z G)$. 

\vp
 
Recently, Juriaans, Passi and Prasad have classified the finite subgroups $G$
whose group $\U_{1}(\Z G)$ is hyperbolic. In the first section we extend
this result, classifying the rings of algebraic integers $R$ of a racional
quadratic extensions and the finite groups $G$ such that
$\U_1(RG)$ is hyperbolic.

\vp

Corrales {\it et all}, in \cite{cjlr}, $2004$, determined
generators of a subgroup of finite index of $\U(\h(\Z(\frac{1+\sqrt{-7}}{2})) )$, whose units have 
norm $1$. 

\vp

For $\h(\Q(\sqrt{d}))$ a division ring we construct some units
of the group $\U(\h (R))$. We obtain a Pell
equation, whose solutions generate the units, which we call
Pell units. Furthermore, we construct units of norm $-1$, which gives
rise to the definition of the Gauss units.


\section{The rings $R$ with $\U_{1}(R G)$ hyperbolic}

Throughout the text, for $d$ a square-free integer we mean that \\ $d\in \D= \{d \in \Z \setminus\{-1,0\}: c^2 \nmid d,$ for all integer $c$ which $c^2 \neq 1\}$. We let $K$ be the quadratic
extension $\Q (\sqrt{-d})$ and $R:=I_K$ be its ring of algebraic integers. The cyclic group of order
$n$ is denoted by $C_n$ and the
quaternion group of order $8$ is denoted $K_8 :=\{\pm 1,\pm i, \pm j, \pm k\} $.

If $G$ is a finite abelian group the unit group $\U_1(RG)$ is a hyperbolic
group if, and only if, its free rank is at most $1$. In \cite{thss}, it is
shown that it is sufficient to consider $G$ a cyclic group of order
$2,3,4,5,6$ or $8$, and  thus  the free rank of $\U_1(RG)$ is calculated. When $G$ is one of the non-abelian  groups of the Theorem $3$ of \cite{jpp}, we show that, in case $\U_1(RG)$ is hyperbolic, $K$ is an imaginary quadratic extension and $G=K_8$. To prove  the converse we use a geometric approach:

\begin{e-definition}Let $K$ be an algebraic number field and $R$ be its ring of algebraic
  integers. For  $a,b \in K$, we denote by  $\h(K)=(\frac{a,b}{K})$ the generalized
  quaternion algebra, i.e.,  $\h(K)$ is the
  $K$-algebra $$\h(K)=K[i,j:i^{2}=a,j^{2}=b,-ji=ij=:k].$$ The set
  $\{1,i,j,k\}$ is a $K$-basis of $\h(K)$. If  $a,b \in R$, then
  $$\h(R)=R[i,j:i^{2}=a,j^{2}=b,-ji=ij=:k].$$  The norm of $x=x_{1}+x_{i}i+x_{j}j+x_{k}k \in \h(K)$ is  $\n(x)=x_{1}^{2}-ax_{i}^{2}-bx_{j}^{2}+abx_{k}^{2}.$ In what follows, we consider $\h(K)=K[i,j:i^{2}=-1,j^{2}=-1,-ji=ij=:k]$.
\end{e-definition}

\begin{e-definition}[\cite{rajw}]
The least natural number $s$ for which the equation
$$-1=a^{2}_{1}+a^{2}_{2}+ \cdots +a^{2}_{s}, a_{j} \in K, 1 \leq j \leq s$$
is soluble is called the stufe of $K$, say $s(K)$. When this equation admits no solution we set $s:=\infty$ and $K$ is called formally real.
\end{e-definition}

Rajwade, in \cite{rajw}, proved that if the quadratic
extension $\Q (\sqrt{-d})$ has $s(K)=4$ then $d \equiv 7 \pmod 8$. Using this, in \cite{thss}, we prove that the quaternion algebra $\h(K)$ over $K$ is a division ring if, and only if, $d \equiv 7 \pmod 8$ and as a corollary we obtain that if $d {\not \equiv} 7 \pmod 8$ then $\U(R K_8)$ is not   hyperbolic.  Defining a proper action of the group $SL_1(\h(R)):=\{x \in \h:\n(x)=1\}$ over the three-dimensional hyperbolic space $\hp$, and a result of Gromov about the fundamental group of a closed $n$-dimension riemannian manifold of constant negative sectional curvature, we prove that if  $d \equiv 7 \pmod 8$ then the group $\U(RK_8)$ is hyperbolic.

\newpage

\begin{theorem}[Theorem $1.7.5$ of \cite{thss}]
Let $R$ be the integral ring of a rational quadratic extension ${K=\Q
 (\sqrt{-d})}$ and $d$ be a square-free integer. The unit group $\U_{1}(R G)$ is
 hyperbolic if, and only if, $G$ is one of the groups listed below
 and $R$ (or $K$) determined by the respective value of $d$:
\begin{enumerate}
 \item $G \in \{ C_{2}, C_{3} \}$ and any $d$.
 \item $G$ is an abelian group of exponent dividing
  $n$ for:\\$n=2$ and $d>0$; or $n=6$ and $d=3$; or $n=4$ and $d=1$.
 \item $G=C_{4}$ and $d>0$.
  \item $G=C_8$ and $d=1$.
  \item $G=K_{8}$ and $s(K)=4$, that is, $d>0$ and ${d \equiv 7 \pmod{8}}$.
\end{enumerate}
\end{theorem}

\vp

 For a
metric space $X$, let the maps $r_1,r_2:[0,\infty[ \longrightarrow
X$ be proper, that is, $r_{i}^{-1}(C)$ is compact for each compact $C \subseteq X$. Two rays are equivalent if for each compact set $C
\subset X$ there exists $N \in \N$, such that, $r_i([N, \infty[),
i=1,2,$ are in the same path connected component of $X \setminus C$. The equivalence
class of $r$ is denoted by $end(r)$; $End(X)$ denotes the set of
equivalence class and $|End(X)|$ is the number of ends of $X$. For a finitely generated group $\Gamma$ and $\G$ its Cayley graph, we
define $Ends(\Gamma):=Ends(\G)$ \cite{grmv}, \cite{mbah}.

\vp

\begin{corollary} \label{qhyp}
 The group $\U(R K_8)$ is hyperbolic if, and
 only if, $d>0$ and  ${d \equiv 7 \pmod
 8}$. Furthermore, the hyperbolic boundary $\partial(\U(R  K_8))\cong S^2$,  the two
 dimensional euclidean sphere, and $\U(R K_8)$ has one end.
\end{corollary}

\vp

Observe that the previous corollary shows a class of hyperbolic groups of one
end which are not virtually free.

\vp

\begin{corollary}\label{frnk}Let $d\equiv 7 \pmod 8$, if $u_1 \cdots u_n \in \U(R K_8)$, then there
  exists $m \in \N$, such that, $\langle u_1^m, \cdots ,u_n^m\rangle$ is a
  free group of rank less ou equal to $n$.
\end{corollary}

\section{The  Pell and  Gauss Units} 
\begin{e-definition}Let $K$ be an algebraic number field and $R$ its ring of algebraic
  integers. For  $a,b \in K$, we denote by  $\h(K)=(\frac{a,b}{K})$ the generalized
  quaternion algebra, i.e.,  $\h(K)$ is the
  $K$-algebra $$\h(K)=K[i,j:i^{2}=a,j^{2}=b,-ji=ij=:k].$$ The set
  $\{1,i,j,k\}$ is a $K$-basis of $\h(K)$. If  $a,b \in R$, then
  $$\h(R)=R[i,j:i^{2}=a,j^{2}=b,-ji=ij=:k].$$  The norm of $x=x_{1}+x_{i}i+x_{j}j+x_{k}k \in \h(K)$ is  $\n(x)=x_{1}^{2}-ax_{i}^{2}-bx_{j}^{2}+abx_{k}^{2}.$ In what follows, we consider $\h(K)=K[i,j:i^{2}=-1,j^{2}=-1,-ji=ij=:k]$.
\end{e-definition}

\vp

\begin{e-proposition} Let $u=u_1+u_ii+u_jj+u_kk \in \U(\h(R))$ with norm $\n(u)$.
 The following conditions hold: 
\begin{enumerate}
\item $u^2=2u_1u-\n(u)$
\item If $d \equiv 7 \pmod 8$ and $\n(u)=1$, then $u$ is torsion if, and only if,
 $u_1\in\{-1,0,1\}$. Thus, the order $o(u)$ is either
 $o(u)=4,2$ or $1$.
\item If $d \equiv 7 \pmod 8$, and $\n(u)=-1$ then $o(u)=\infty$. 
\end{enumerate}
\end{e-proposition}

\vp

Let $\q:=\Q (\sqrt{d})$  and $\xi \neq \psi \in
 \{1,i,j,k\}$. For $\epsilon=x +y \sqrt{d} \in \U(I_\q)$, we denote $u_{(\epsilon)}:=x \sqrt{-d}\xi + y \psi \in \h(K)$.
\vp 

\begin{e-proposition}
Let $d \equiv i \pmod 4, i \in \{2,3\}$ and $\xi \neq \psi \in
 \{1,i,j,k\}$. The following conditions hold:
\begin{enumerate}
\item  $u_{(\epsilon)} \in \U(\h(R))$ if, and only if, $\epsilon =p + m\sqrt{d} \in \U(I_{\q})$. 
\item If $1 \notin supp(u)$ then $u_{(\epsilon)}$ is torsion.

\item If $\mu,\nu \in \U(I_\q)$ and $1 \in supp(u_{(\mu)}) \cap supp(u_{(\nu)})$, then $u_{(\mu)} u_{(\nu)}=u_{(\mu \nu)}$.
\item If $1 \in supp(u_{(\epsilon)})$, then $\langle u_{(\epsilon)} \rangle = \{u_{(\epsilon^n)}, n \in \Z\}$.
\item For $d \equiv 3 \pmod 4$ and $\f:=\Q(\sqrt{2d})$.  $$u=m\sqrt{-d} \xi + p \psi + (1-p) \phi \in \U(\h(R))
\Leftrightarrow \epsilon=(2p-1)+m\sqrt{2d} \in \U(I_{\f})$$
\end{enumerate}
\end{e-proposition}

\vp

\begin{theorem} Let $\h(K)$ be a division ring. If
  $x+y\sqrt{d} \in \U(I_\q)$, then $$u=\left \{ \begin{array}{ll} \frac{y}{2}\sqrt{-d}+(\frac{y}{2}\sqrt{-d})i+(\frac{1 \pm x}{2})j+(\frac{1 \mp x}{2})k& \textrm{if } y \equiv 0 \pmod 2\\
  xy\sqrt{-d}+(xy\sqrt{-d})i+(\frac{1 \pm (x^{2}-y^{2}d)}{2})j+(\frac{1 \mp
  (x^{2}-y^{2}d)}{2}k)& \textrm{if } y \equiv 1 \pmod 2  \end{array} \right.$$
are units in $\h(R)$.
\end{theorem}

\vp

\begin{e-definition}
The given units above are called \textbf{Pell Units}. For $l\in\{2,3\}$, a
\textbf{Pell $l$-unit} is a unit whose support has cardinaliy $l$, and the
unique non integer coefficient is of the form $m\sqrt{-d}$. 
\end{e-definition}

\vp
\begin{theorem}
Let  $\h(K)$ be a division ring. If $m \equiv 2 \pmod 4$, then there exist integers
$p,q,r$, such that, $u=m \sqrt{-d}+pi+qj+rk \in \U(\h(R))$.
\end{theorem}

\vp

\begin{e-definition}
A unit $u$  of $\h(R)$ whose support has cardinality $l:=|supp(u)|>1$, the unique
non integer coefficient of $u$ is of the form  $m\sqrt{-d}$ and $m^2 \pm 1$ is a sum of three square integers is called a \textbf{Gauss unit}, or a Gauss $l$-unit.
\end{e-definition}

\vp

\begin{e-proposition}Let $u$ be a unit of norm $\n(u)=1$, $l\in \{2,3\}$, and
 $\h(K)$ a division ring.
$$\textrm{$u$ is a Pell $l$-unit if, and only if, $u$ is a Gauss $l$-unit}$$
\end{e-proposition}

\vp

\begin{theorem}Let $d \equiv 7 \pmod 8$. If  $u,v \in
\U(\h(R))$ are Gauss $2$-units, and $supp(u)\cap supp(v)=\{1\}$, then there exists $m \in \N$, such that, $\langle u^m,v^m \rangle$ is a free group of rank two. 
\end{theorem}

\vp

In \cite{cjlr}, the authors exhibit a set of generators $S$ of $SL_1(\h(\Z (\frac{1+\sqrt{-7}}{2})))$. The gauss unit $v=6\sqrt{-7}+15i+5j+k$ has norm $\n(v)=-1$, therefore $\U(\h((\frac{1+\sqrt{-7}}{2})))=\langle S,v \rangle$.
The elements of the set $S\setminus \{i,j\}$ are units of the form $\frac{m +
  \sqrt{-7}}{2}+(\frac{m - \sqrt{-7}}{2})i+pj$, possibly with a permutation of
the coefficients. If the condition $d \equiv 7 \pmod 8$ is assumed. Then the
solutions of the equation $m^2+2p^2=2+d$, give rise to these units.

\vp

\section*{Acknowledgements}
This work is part of the second authors Ph.D
thesis. He would like to thank his thesis supervisor Prof. Dr. Stanley
Orlando Juriaans for his guidance during this work.

\end{document}